\documentclass[english]{lni}

\usepackage{booktabs}

\usepackage[]{blindtext}

\usepackage[english]{babel}
\DefineBibliographyStrings{english}{
  references   = {References},
  bibliography = {References},
}

\begin{document}

\title[K.I.A Derouiche]{Zygmunt Janiszewski’s Mathematical Publications During His Engagement in Poland’s Struggle for National Independence in World War I (Les publications mathématiques de Zygmunt Janiszewski pendant son engagement pour l’indépendance nationale en Pologne pendant la Première Guerre mondiale)}
\author[1]{K.I.A Derouiche}{kamel.derouiche@gmail.com,k.derouiche@algeriemotors.com}{0000-0002-2311-5113}
\affil[1]{BMW Motorrad Algeria – IT Department, Rue Mohammed Loubi 08
		16040 Hussein Dey, Alger, Algérie}
\maketitle

\begin{abstract}
This article focuses on the mathematical publications of Zygmunt Janiszewski (1888–1920), a major figure in Polish science at the beginning of the twentieth century. Serving in the Polish Legion between 1914 and 1920 in the struggle for national independence, Janiszewski was not only one of the founders of the Polish School of Mathematics, but also the initiator of the creation of Fundamenta Mathematicae, the world’s first specialized scientific journal devoted exclusively to pure mathematics.

The aim of this study is to reconstruct, as far as the available sources allow, the circumstances and conditions under which Janiszewski published his mathematical works. Through this analysis, we highlight the significance of his research in the broader process of structuring and consolidating the modern Polish mathematical community.

	\selectlanguage{Polish}
\textbf{Abstrakt} Artykuł poświęcony jest publikacjom matematycznym Zygmunta Janiszewskiego 
(1888–1920), jednej z najwybitniejszych postaci polskiej nauki początku XX wieku. W latach 
1914–1920 Janiszewski służył w Legionach Polskich, walcząc o niepodległość kraju. Był nie tylko 
współtwórcą polskiej szkoły matematycznej, lecz także inicjatorem powstania czasopisma 
Fundamenta Mathematicae — pierwszego na świecie specjalistycznego pisma naukowego poświęconego matematyce czystej.
\\

Celem niniejszego opracowania jest zbadanie — w miarę możliwości źródłowych — okoliczności i 
warunków, w jakich Janiszewski publikował swoje prace matematyczne. Analiza ta ukazuje 
znaczenie jego badań w procesie kształtowania i konsolidacji nowoczesnej polskiej społeczności 
matematycznej.

	\selectlanguage{french}
\textbf{Résumé}. Cet article est consacré aux publications mathématiques de Zygmunt Janiszewski (1888–1920), figure majeure de la science polonaise du début du XXeme siècle. Engagé dans la Légion polonaise pour l’indépendance nationale entre 1914 et 1920, Janiszewski fut non seulement l’un des fondateurs de l’école polonaise de mathématiques, mais également l’initiateur de la création de la revue Fundamenta Mathematicae, première revue scientifique spécialisée au monde consacrée aux mathématiques pures.

L’objectif de cette étude est d’examiner, autant que les sources le permettent, les circonstances et les conditions dans lesquelles Janiszewski publia ses travaux mathématiques. À travers cette analyse, nous mettons en lumière la portée de ses recherches dans le processus de structuration et de consolidation de la communauté mathématique polonaise moderne
\end{abstract}


\begin{center}
	\textbf{{\large Preprint (Version Courte)}}
\end{center}

\section{Introduction}
The purpose of this article is to examine the circumstances and historical context in which the 
mathematician, philosopher, and co-founder of the modern Polish school of mathematics, Zygmunt 
Janiszewski, conducted and published his research during the Great War (1914–1917), at a time 
when he was simultaneously serving as a legionnaire in the Polish Legions.
\par

Celem niniejszego artykułu jest zbadanie okoliczności oraz kontekstu historycznego, w których matematyk, filozof i współtwórca nowoczesnej polskiej szkoły matematycznej, Zygmunt Janiszewski, prowadził i publikował swoje badania podczas Wielkiej Wojny (1914–1917), pełniąc równocześnie służbę w Legionach Polskich.
\par

\textbf{Le contexte historique}. Lorsque la Première Guerre Mondiale
éclate le 28 juillet 1914, elle oppose l'Allemagne et l'Autriche-Hongrie
à la Russie soutenue par la France et la Grande-Bretagne. Les polonais
pouvaient soit rester passifs, soit prendre parti dans le conflit.
Beaucoup furent enrôlés dans leurs armées respectives certains ne
croyaient pas à mener une guerre servant les intérêts des puissances
occupantes ; une guerre considérée comme fratricide, compte tenu de la
présence polonaise dans les États ennemis. D'un autre côté, beaucoup ont
vu une opportunité de se battre pour l'état polonais indépendant en
s'alignant du côté des puissances centrales ou de la Triple Entente. Cela
était particulièrement vrai en Galicie et pour d'autres régions,  où de
jeunes hommes se sont massivement portés volontaires ont été enrôlés par
forces dans les légions polonaises parmi eux des mathématiciens, qui ont
servi dans les tranchées ou dans des unités non combattantes, subissant
des blessures, des empoisonnements au gaz, des camps de prisonniers ou
l'internement. Ils ont également été touchés par les évacuations forcées,
les restrictions de voyage, les pénuries de nourriture et de matières
premières. Aucun n'a perdu la vie à la suite d'opérations de guerre, mais
deux universitaires exceptionnels le physicien Marian
Smoluchowski(1872-1917) et Zygmunt Janiszewski sont morts d'épidémies du
à la grippe espagnole deux années après la fin de première mondiale. Bien
que les destins que nous décrivons aient été typiques à bien des égards,
présenter les réalités de la guerre d'un point de vue personnel n'est pas
le seul but de notre article. Nous voulons plutôt nous concentrer sur
l'impact de guerre et les conditions qu’elle engendres sur les activités
de recherches pour les mathématiciens. Notons que certains mathématiciens
étendent leur activité de guerre à l'enseignement au secondaire ou au
primaire quand et où le besoin s'en fait sentir. Quelques-uns se sont
également engagés en dehors des mathématiques et de l'éducation : dans
des activités politiques Wiktor Staniewicz(1866-1932), dans l'expression
artistique Leon Chwistek(1884-1944) ou dans l'écriture sur des thèmes
culturels et religieux juifs Chaim Müntz(1884-1956). Outre les chercheurs
connus pour leurs résultats exceptionnels en mathématiques avant la
guerre ou après, tels que Banach et Sierpiński, nous présentons de
nombreuses personnes qui ont apporté des contributions moindres aux
connaissances mathématiques, en particulier celles (par exemple Izabela
Abramowicz, Zygmunt Chwiałkowski, Adam Patryn) qui n'ont pas poursuivi
leurs recherches après la guerre. Enfin, nous parlons de certains
physiciens, astronomes, ingénieurs et philosophes qui, dans les
circonstances de la guerre, se sont engagés dans l'enseignement des
mathématiques au niveau académique ou dans les activités des sociétés
savantes aux côtés de leurs collègues mathématiciens. Bien sûr, nous ne
pouvons pas citer tout le monde, surtout si les informations concernant
leurs années de guerre sont rares. Nous avons même omis certains
mathématiciens bien connus Franciszek Mertens, Stefan Bergman et autres
pour qui la guerre s'est produite pendant leurs années d'âge adulte ou
senior. Dans l'ensemble, les cours académiques se sont poursuivis avec
des interruptions inévitables même si de nombreux étudiants et
professeurs ont servi dans l'armée, des bâtiments ont été réquisitionnés
à des fins et ressources militaires (bibliothèques, matériel
scientifique, etc.) ont été évacués. Des doctorats ont été décernés
Franciszek Leja, Witold Wilkosz, Antoni Plamitzer, Adam Patryn ou reconnu
Zygmunt Janiszewski; habilitations ont été accordées Hugo Steinhaus,
Stanisław Ruziewicz, Antoni Łomnicki, Eustachy Żyliński, Tadeusz
Banachiewicz, Stanisław Leśniewski ou nié Lucjan Emil Bottcher. Beaucoup
de nouvelles personnes talentueuses ont été attirées par les
mathématiques. Un cas remarquable est celui de Stefan Banach,  ayant une
carrière mathématique a été stimulée par sa rencontre fortuite avec Hugo
Steinhaus en 1916, ce qui n'arriverait probablement pas si Steinhaus
n'occupait pas un poste administratif à Cracovie après sa libération de
l'armée\citet{Domoradzki2018}. Revenons au cas exceptionnels
de Zygmunt Janiszewski l'un des fondateurs de l'école mathématique
polonaise a émergé et a gagné le monde entier reconnaissance. Il convient
de s’interroger sur la manière dont Zygmunt Janiszewski parvint à
réaliser une œuvre scientifique d’une telle envergure dans un contexte de
guerre particulièrement contraignant et en un temps remarquablement
limité. On peut également se demander dans quelle mesure il aurait pu
poursuivre sa carrière universitaire après le conflit, compte tenu des
graves difficultés matérielles et institutionnelles engendrées par cette
période.
\subsection{Les sources}
À ce stade, le papier est en cours de rédaction. Nous nous appuyons sur
diverses sources : outre les articles publiés par Zygmunt Janiszewski
lui-même, les travaux des mathématicien(ne)s polonais Bronisław Knaster,
Stanisław Domoradzki et Małgorzata Przeniosło, ainsi que l’ensemble des
dossiers biographiques, photographies et documents militaires conservés
aux archives du bureau d’histoire militaire polonais.

\subsection{Organisation de l'article}
La section 2 sera consacrée à la biographie de Zygmunt Janiszewski ; la
section 3 traitera de sa scolarité et de ses études universitaires ; la
section 4 exposera la chronologie de sa mobilisation dans la légion
polonaise ; la section~5 abordera les différents articles scientifiques
publiés durant la période de sa mobilisation. Enfin, l’article se
conclura dans la section 5 par l’analyse des circonstances de la fin de
sa mobilisation dans la légion polonaise, dans la section 6 l'article se concluera .

\section{Quelques éléments biographiques sur Zygmunt Janiszewski}
Zygmunt Janiszewski est né le 12 juin 1888 à Lwów (Lemberg en allemand,
Lviv en ukrainien) \cite{Bechtel2005}. Il était enfant unique d’une
famille de la petite bourgeoisie. Son père, Czesław Janiszewski, était
diplômé de l’École d’économie de Varsovie, fonctionnaire au bureau du
procureur général du Royaume de Pologne, puis directeur de la banque
Towarzystwo Kredytowe Miejskie (la Société de Crédit Municipal) à
Varsovie; sa mère, Julia Szulc-Chojnicka, fait l’objet de peu
d’informations biographiques ; on sait toutefois d’après les informations
que nous est parvenu que ses parents étaient décédés avant 1916
\cite{Przeniosło2011}.

Dès janvier 1905, immédiatement après sa participation à la grève
scolaire \cite{Zarnowska1986}, le jeune Zygmunt Janiszewski s’installe à
Lwów, où il passe son diplôme d’études secondaires.
En 1907, il obtient ce diplôme et décide de se rendre à Zurich à
l’automne de la même année pour entreprendre des études à la faculté de
philosophie, comme le mentionne son collègue de l’époque, le
mathématicien Stefan Straszewicz (1898–1927) \cite{Przeniosło2011}.
Il déclarera plus tard qu’il n’y trouva pas des conditions favorables au
développement de la créativité et de la pensée mathématique indépendante.

\section{Scolarité et études universitaires}
Dès janvier 1905, immédiatement après sa participation à la grève
scolaire sur le territoire du Royaume de Pologne \cite{Ger1986}, le jeune
Zygmunt Janiszewski s’installe à Lwów, où il passe son diplôme d’études
secondaires. 
En 1907, après l’obtention de ce diplôme, il part pour Zurich à l’automne
afin d’entamer des études à la faculté de philosophie, comme le mentionne
Stefan Straszewicz \citet{Przeniosło2011}.
\\
Par la suite, Janiszewski est transféré au deuxième semestre de la
première année d’études à l’université de Göttingen, où il rencontre D.
Hilbert, Minkowski et Zermelo.
\par

Après cette période de formation à Göttingen \citet{Przeniosło2011},
Janiszewski s’installe à Paris pour soutenir sa thèse de doctorat en
topologie, «\,Sur les continus irréductibles entre deux points\,»
\citet{Janiszewski1911}, devant un jury composé d’éminents mathématiciens
français, dont Henri Poincaré et Maurice Fréchet, avec lesquels il
conserva une amitié.

\section{Mobilisation — Zygmunt Janiszewski à la Légion polonaise}
Après le déclenchement de la Première Guerre mondiale, Zygmunt
Janiszewski rejoint la Légion de l’Est formée à Lwów le 8 août ; celle-ci
est dissoute le 21 septembre.
Il servit comme artilleur dans le 1\textsuperscript{er} escadron
d’artillerie commandé par Mieczysław Gałusiński \footnote{Mieczysław
Gałusiński (d\'ec\'ed\'e après 1935), pseudonyme Jełowiecki, \'ecrivain,
journaliste, co-organisateur de l'artillerie l\'egionnaire, il est promu au
grade de capitaine dans les Légions.}.

L’escadron se composait alors de trois batteries et n’était pas bien
armé. Il rejoint les formations de légion commandées par le général Karol
Trzaska-Durski (1849–1935), envoyées en première ligne dans les Carpates
orientales. Le transport ferroviaire partit de Cracovie pour Huszt (Hongrie) le 
1\textsuperscript{er} octobre 1914 \cite{Przeniosło2011}.
\begin{figure}[!tb]
	\centering
		\includegraphics[width=0.6\textwidth]{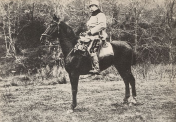}
		\caption{Le générale Karol Durski-Trzaska\cite{Legionów1916} }
		\label{fig:KarolDurski-Trzaska}
\end{figure}

Quelques unités livrèrent les premières batailles quelques jours plus
tard dans la région de Huszt et Marmaros-Sziget, puis la plupart des
formations légionnaires se dirigèrent vers la frontière galicienne avec
l’intention de capturer les environs de Rafajłowa, Zielona et Nadwórna
(près de Stanisławów).
Parmi elles se trouvait la 2\textsuperscript{e} batterie d’artillerie,
tandis que la 1\textsuperscript{re} batterie fut mise à la disposition de
la 52\textsuperscript{e} division d’infanterie autrichienne.
\par
Il est difficile de préciser dans laquelle de ces deux batteries
Janiszewski a servi. Il existe un document confirmant l’affectation de
Janiszewski à deux batteries, délivré le 16 décembre 1914 par le
commandant de toute l’escadre, le prêtre Jełowiecki, tandis que Wacław
Chocianowicz le mentionne comme faisant partie de la
1\textsuperscript{re} batterie.
Il est possible que, pendant le service de Janiszewski, l’affectation ait
été modifiée ; par conséquent, des informations sur les deux unités sont
brièvement présentées ci-dessous. 
\par
La 2\textsuperscript{e} batterie d’artillerie, avec d’autres formations
légionnaires, était en Galicie le 21 octobre. Le franchissement de la
frontière impliquait le passage des cols des Carpates.
Voulant surprendre les Russes, le commandement des Légions prit la
décision audacieuse de prolonger le chemin étroit de cinq kilomètres
menant à travers le col de Rogoga, au pied du mont Pantyr ; la route fut
construite en rondins posés sur une pente préalablement profilée.

Après avoir traversé la frontière, les troupes légionnaires capturèrent
Rafajłowa et Zielona, puis livrèrent de féroces batailles pour Nadwórna
et ses environs. 
La 2\textsuperscript{e} batterie combattit notamment à Hwozd le 26
octobre et trois jours plus tard à Mołotków.
\par
Après ces pertes significatives durant la bataille, les troupes qui y
participèrent furent retirées en Hongrie pour un court repos, puis
réoccupèrent les positions déjà acquises. À cette époque, la
1\textsuperscript{re} batterie participait aux opérations autrichiennes
dans les environs de Mikuliczyn.
\par
Auparavant, elle avait effectué près de deux semaines de marches
épuisantes (un total d’environ 200 km) vers diverses positions en
constante évolution ; elle ne rejoignit les formations de la Légion mère
que le 16 novembre. Fin novembre, les unités légionnaires furent divisées
en deux groupes tactiques. 
\begin{figure}[!tb]
	\centering
		\includegraphics[width=0.4\textwidth]{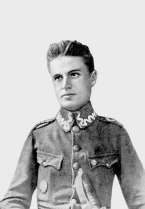}
		\caption{La seule image qui nous est parvenue de Janiszewski en tenue militaire }
		\label{fig:KarolDurski-Trzaska2}
\end{figure}

\section{Publications de Zygmunt Janiszewski durant la période 1914–1917}
Durant la période de la Première Guerre mondiale, Zygmunt Janiszewski parvint à publier
plusieurs articles importants malgré les conditions difficiles liées au conflit.
Ses publications reflètent une double orientation : d’une part, ses travaux strictement 
mathématiques ; d’autre part, ses réflexions philosophiques et pédagogiques sur la science en 
Pologne.

\begin{table}[h!]
	\centering
	\caption{Publications de Zygmunt Janiszewski entre 1914 et 1917}
	\label{tab:publications}
	\begin{tabular}{|l|l|l|}
		\hline
		\textbf{Année} & \textbf{Titre / Publication} & \textbf{Type de contribution} \\ \hline
		1914 & O realizmie i idealizmie w matematyce & Article philosophique \\ 
		1915 & Contribution au \textit{Poradnik dla samouków} & Chapitre d’ouvrage \\
		1916 & Notes pédagogiques diverses (Varsovie) & Articles courts \\ 
		1917 & Propositions pour une revue de mathématiques polonaise & Texte programmatique \\ 
		\hline
	\end{tabular}
\end{table}

\subsection{La contribution au \textit{Poradnik dla samouków} (« Guide pour autodidactes »)}
L’ouvrage collectif \textit{Poradnik dla samouków} (« Guide pour
autodidactes »), dirigé par Aleksander Heflich, fut publié en 1915 à
Varsovie, c’est-à-dire avant la mobilisation complète de Janiszewski dans
la Légion polonaise.
Janiszewski y contribua par un texte remarquable sur la structure de la
recherche mathématique et la place de la logique dans la formation de
l’esprit scientifique.
Cet ouvrage avait pour but d’offrir aux autodidactes polonais un panorama
méthodologique de la culture scientifique, à un moment où l’accès aux
études supérieures était limité par les conditions politiques.

La contribution de Janiszewski s’inscrit pleinement dans sa conception
d’une science organisée, collective et nationale ; elle préfigure
également ses appels ultérieurs à la création d’un centre de recherche et
d’une revue polonaise spécialisée en mathématiques.

Le 30 août 1914, Janiszewski s'enrôle dans les légions polonaises. Il
participe à la campagne des Carpates en 1914/1915. Lorsque les Allemands
prirent Varsovie, Janiszewski y fut appelé à un poste. Il raconta plus
tard à Hugo Steinhaus que, lorsqu’il arrivait, son chauffeur de la gare
au commandement était le géomètre Max Dehn (1878–1952), ancien élève de
Hilbert à Göttingen, Le service militaire de Janiszewski (alors déjà à
portée de un sergent). A cette période et l'internement de Wacław
Le fasicule « Poradnik samoukków »  (« Conseiller pour les
autodidactes »)\cite{Janiszewski1915b} est une série d’articles sur la
philosophie et les différentes branches des mathématiques destinée aux
autodidactes, et à ceux qui voudraient accompagner ou organiser un
apprentissage autonome. Le premier tome dédié à un état de l’art aux
diverses branches des mathématiques de son époque
, des introductions méthodologiques, des indications bibliographiques, et
des commentaires pédagogiques sur les niveaux d’étude qui témoigne de
l’énorme influence exercé par Janiszewski sur le développement des
mathématiques en Pologne, ce guide est structuré en plusieurs niveaux
d’étude : élémentaire, moyen, et « IIIeme degré » . Ces niveaux servent à
organiser la progression de l’apprenant autodidacte, chaque section est
présenté  comme une sous-discipline mathématique qui commence par des
introductions méthodologiques, des définitions, un exposé des contenus,
références, éventuellement des exercices ou indications bibliographiques
critiques. Notre mathématicien est l’auteur de plusieurs chapitres clés :
introduction générale aux mathématiques, fondements de la géométrie,
topologie, équations différentielles et fonctionnelles, séries,
philosophie des mathématiques, conclusion, et une section sur les
possibilités d’études à l’étranger et des données sur les universités
européennes. Historiquement Poradnik dla samouków a joué un rôle clé dans
l’éducation scientifique autodidacte en Pologne avant que les
institutions universitaires modernes. Il servait à combler les lacunes
d’accès, notamment dans les régions éloignés ou par rapport aux
publications étrangères.  Pour Janiszewski, ce guide a été aussi un moyen
de diffuser ses idées fondamentalistes et topologiques à un public plus
large, pas seulement aux spécialistes.
\begin{figure}[!tb]
	\centering
		\includegraphics[width=0.4\textwidth]{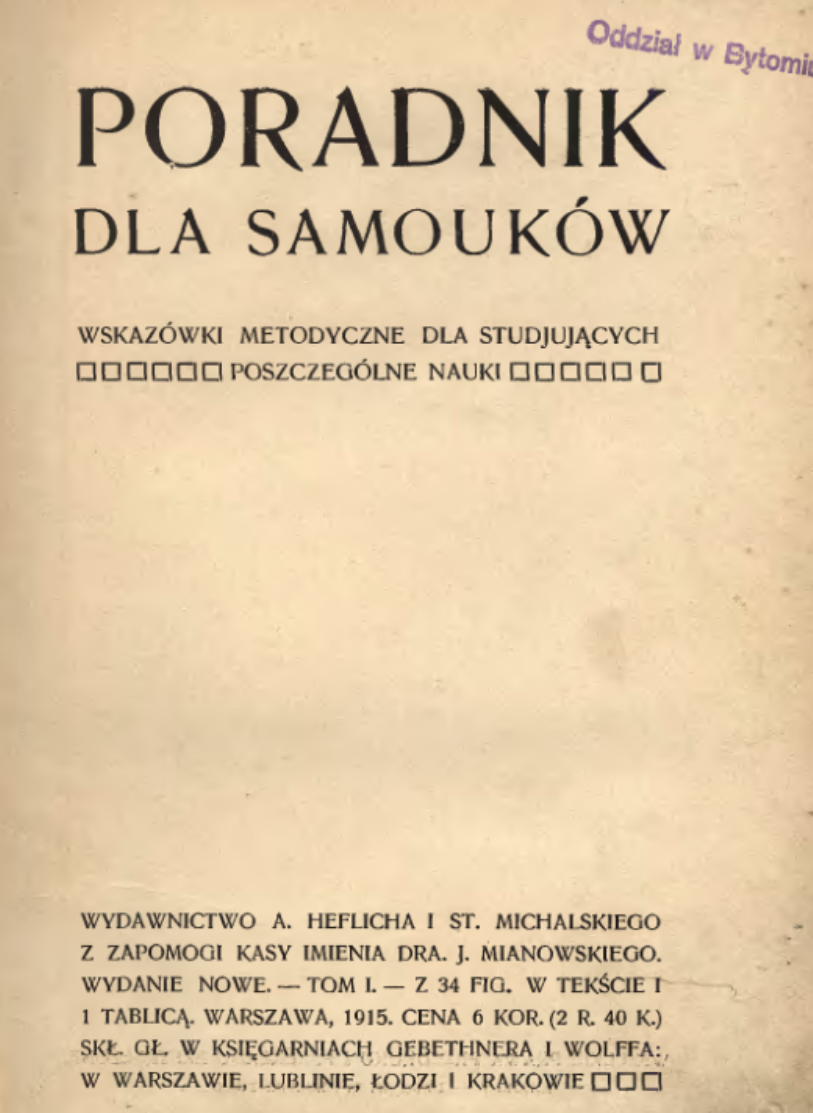}
		\caption{Première page du manuel consacré aux autodidactes}
		\label{fig:Premire page de l'article}
\end{figure}

Poradnik dla samouków\cite{Janiszewski1915b} a joué un rôle clé dans
l’éducation scientifique autodidacte en Pologne avant que les
institutions universitaires modernes ne se soient pleinement
stabilisées. Il servait à combler les lacunes d’accès, notamment dans les
régions ou par rapport aux publications étrangères, pour Janiszewski, ce
guide a été aussi un moyen de diffuser ses idées fondamentalistes et
topologiques à un public plus large, pas seulement spécialisé.
\par

\begin{table}[h!]
	\centering
	\caption{Contributions de Zygmunt Janiszewski dans le guide Poradnik
	dla samouków}
		\label{tab:contributions}
\begin{tabular}{|p{5cm}|p{10cm}|}
	\hline
	Rubrique & Contribution apporté par Janiszewski \\
	\hline
  Introduction générale  & Mise en perspective de ce qu’est la discipline
  mathématique, son rôle parmi les sciences, importance de la
  méthodologie, guidance pour l’autodidacte. \\
    \hline
  Introduction au III \^eme degré & Suggestions pour les apprenants
  avancés : quels domaines explorer, comment structurer ses études pour
  atteindre un niveau universitaire \\
	\hline
	Développements : Équations différentielles ordinaires, équations
	fonctionnelles/différentielles/intégrales, séries
	& Présentation claire des thèmes, des définitions, des méthodes de
	résolution, des résultats classiques jusqu’alors, avec des références.
	Cela montre son effort de synthèse dans des domaines techniques.
	\\
	\hline
	Topologie & Chapitre sur la topologie, exposant des notions comme
	connexité, ensembles, séparation, etc., permettant au lecteur
	autodidacte de comprendre les concepts modernes de la topologie
	naissante
	\\
	\hline
	Fondements de la géométrie & Étude des axiomes géométriques, de la
	logique géométrique, de la justiciabilité des définitions, ce qui
	reflète les préoccupations fondationnalistes de Janiszewski. \\
	\hline
	Logique/Philosophie des mathématiques & Réflexions sur la nature des
	mathématiques, leur méthode, leur philosophie, les principes sous-
	tendant les définitions et axiomes. \\
	\hline
	Conclusion & Synthèse des connaissances exposées, conseils pour le
	travail scientifique, des renseignements sur les universités
	européennes, ce qui montre non seulement la connaissance de la
	littérature internationale mais aussi une vision éducative large.\\
	\hline
\end{tabular}
\end{table}

Au cours de l'année universitaire 1915/16. Cela a rendu très difficile
pour l'université d'offrir des conférences et des séminaires réguliers en
mathématiques. En raison de cette situation difficile, en octobre 1916,
la Faculté de philosophie demanda au recteur de demander le rappel de
Janiszewski

Les conditions de publication de ce travail furent particulièrement
difficiles : les interruptions postales, la censure et les déplacements
constants rendaient l’échange de manuscrits laborieux.
Janiszewski bénéficia toutefois du soutien de ses collègues de Varsovie
et de Lwów, notamment Wacław Sierpiński et Stefan Mazurkiewicz, avec
lesquels il entretenait une correspondance scientifique intense.

\subsection{Vers la fondation de \textit{Fundamenta Mathematicae}}
Durant les années de guerre, Janiszewski réfléchit à l’organisation de la
recherche mathématique en Pologne. Dans son texte programmatique de 1917,
il proposa la création d’une revue scientifique spécialisée, destinée à
rassembler les travaux de l’école polonaise renaissante.
Cette idée aboutit, après sa mort, à la fondation en 1920 de la revue
\textit{Fundamenta Mathematicae}, première publication au monde consacrée
exclusivement aux mathématiques pures.
Son collègue Stefan Mazurkiewicz, avec Wacław Sierpiński, prit en charge
la réalisation de ce projet.
Ainsi, Janiszewski, bien que disparu prématurément, demeure le véritable
initiateur de cette entreprise scientifique collective.

\subsection{Son article sur la philosophie mathématique :  O realizmie i idealizmie w matematyce }
L’année 1916, est une période troublée au plein Première Guerre Mondiale,
les mathématiciens polonais commençaient à organiser ce qui sera l’«
école polonaise » (programme janiszewskien : concentration sur théorie
des ensembles, topologie, logique).  Le texte "O realizmie i idealizmie w
matematyce "(Sur le réalisme et l'idéalisme en mathématiques)
\citet{Janiszewski1916} s’inscrit donc à la croisée de préoccupations
purement philosophiques (ontologie des entités mathématiques) et de
questions très concrètes de pratique mathématique émergente (p.ex.
axiomatique, choix d’axiomes) ,  Janiszewski pose la question qui
tourment philosophes et mathématiciens depuis Platon : les objets
mathématiques existent-ils « là-hors » indépendamment de l’esprit
(réalisme/platonisme) ou sont-ils constructions mentales/formelles
(idéalisme / néo-kantisme / constructivisme) ?\cite{Murawski2004};
Plutôt que d’opter pour un dogme, il cherche à articuler une position
méthodologique qui mette en valeur la pratique mathématique (les méthodes
acceptées, l’usage d’axiomes) tout en reconnaissant la force des
intuitions  dite « objectivantes » qui guident les mathématiciens quand
pourrait résumé en les points suivants:
\begin{enumerate}
	\item Les entités mathématiques ne sont pas des “objets matériels” mais
	possèdent une réalité opératoire : elles acquièrent statut ontologique
	par leur rôle stable et reproductible dans des systèmes déductifs
	cohérents.
	\item Autrement dit, l’« existence » en mathématiques est conditionnée
	par la fonction et la cohésion structurelle dans un réseau d’énoncés
	(et non par une hypothétique contact métaphysique avec un domaine
	platonesque). Cette position est à la fois critique envers un réalisme
	naïf et distante d’un idéalisme pur.
\end{enumerate}
Arguments et éléments du raisonnement (exposés selon Janiszewski).
\begin{enumerate}
	\item Différence entre existence ontologique et existence opératoire,
	Janiszewski insiste sur la séparation conceptuelle entre « exister »
	métaphysiquement et « être utile/nécessaire » au sein d’un système
	mathématique. La valeur d’un concept se mesure par sa place et sa
	nécessité dans les inférences.
	\item Rôle des axiomes et de la pratique, il montre que le choix ou le
	refus d’un axiome (par exemple, l’acceptation de certains principes non
	constructifs) est à la fois philosophique et méthodologique : des
	écoles « réalistes » et « idéalistes » peuvent diverger sur des axiomes
	sans que l’unique critère soit métaphysique. Sur ce point historique,
	Janiszewski évoque explicitement les débats autour d’axiomes (notamment
	le rôle que jouera plus tard le débat sur l’axiome du choix).
	\item Cohérence et nécessité logique comme critère de vérité
	mathématique,  loin d’une simple subjectivité, la vérité mathématique
	se rattache à la cohérence et à la déductibilité interne ; c’est cette
	« nécessité formelle » qui explique pourquoi les mathématiques donnent
	l’impression d’objectivité.
\end{enumerate}
Observations techniques/points précis (ce qui fait encore écho
aujourd’hui).

\begin{enumerate}
\item  Axiomatique vs intuition : Janiszewski note qu’une partie du
désaccord philosophique tient aux préférences méthodologiques
(acceptation d’axiomes non constructifs, usage intensif du raisonnement
transﬁnit), pas seulement à un choix ontologique pur. Ceci anticipe les
débats modernes sur le statut d’axiomes comme l’axiome du choix ou
certaines formes du principe du tiers exclu.
	\item Éléments structuraux : sa description de la « réalité des entités
	» comme dépendant de leur rôle dans un système mathématique se lit
	aujourd’hui comme une préfiguration des thèmes structuralistes (la
	primauté des structures et des relations sur la nature intrinsèque des
	objets). Janiszewski ne formule pas la théorie structuraliste moderne,
	mais ses intuitions méthodologiques lui sont voisines.
\end{enumerate}
L'article à eu un impact sur la jeune école polonaise dans son
développement et évolution par des noms tel-que Kuratowski (1896-1980),
Mazurkiewicz(1888-1945), Leśniewski (1886–1939) , Alfred Tarski
(1901-1983) et Kotarbiński (1886 - 1981) à la philosophie des
mathématiques : il alimenta la culture intellectuelle qui rendit
possible, dans les années 1920–30. Janiszewski, par son programme et ses
prises de position, aida à définir les orientations (fondements, logique,
topologie) privilégiées ensuite par la communauté
polonaise\cite{Woleński2022}, ce texte témoigne de la réflexion
épistémologique profonde de Janiszewski sur la nature ontologique des
objets mathématiques et sur les fondements du savoir mathématique.

\begin{figure}[!tb]
	\centering
		\includegraphics[width=0.6\textwidth]{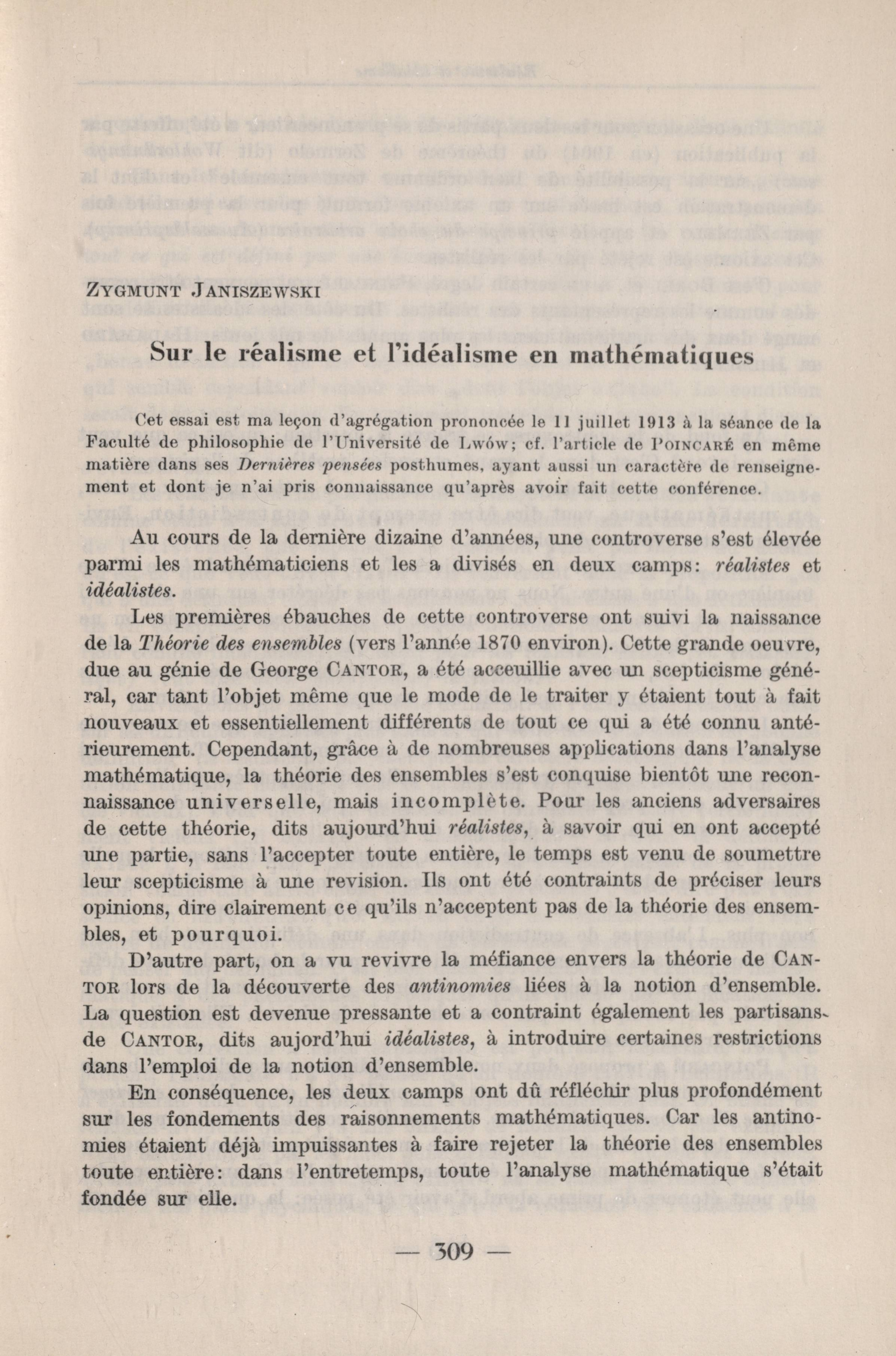}
		\caption{Article de Janiszewski sur la philosophie des mathématiques}
		\label{fig:Premire page de l'article}
\end{figure}

\section{La fin de sa mobilisation dans la légion polonaise}
En juillet 1917, après avoir refusé de prêter le serment qui obligeait
les légionnaires à défendre fidèlement la fraternité d’armes avec les
troupes d’Allemagne et d’Autriche-Hongrie et avec leurs pays alliés,
cette décision eut de lourdes conséquences.
À la suite de la « crise du serment » et de la dissolution des Légions
polonaises\cite{Norman1982}, Janiszewski dut quitter l’université, bien
que sa nomination en tant que maître de conférences adjoint eût été
prolongée de deux ans, du 1\textsuperscript{er} octobre 1917 au 30
septembre 1919.
Durant cette période, il se cachait dans le Royaume de Pologne,
principalement à Radom et Kielce (Galicie), sous le nom d’emprunt de
Zygmunt Wicherkiewicz.
Il séjourna d’abord dans le village de Boiska, dans le poviat d’Iłża, où
il s’impliqua dans le soutien organisationnel des écoles locales. %
corrigé : virgules et accords
Il trouva cet endroit probablement grâce aux connaissances acquises lors
de ses activités dans le bureau central de l’école à Piotrków.
De nombreuses informations sur cette période et les années suivantes sont
fournies par ses lettres à sa fiancée Janina Kelles-Krauz (1898–1975),
qu’il connaissait depuis son séjour à Sports Fields.
Zygmunt Janiszewski est mort subitement à Lviv le 3 janvier 1920, à l’âge
de 31 ans, de la grippe espagnole qui toucha la Pologne entre 1918 et
1920.
\section{Conclusion}
L’examen du parcours scientifique et militaire de Zygmunt Janiszewski met
en lumière le destin singulier d’un savant profondément engagé dans la
vie intellectuelle et nationale de son pays.
Malgré les épreuves de la guerre, les privations et la maladie, il
contribua à poser les bases de la topologie moderne et de la philosophie
mathématique en Pologne.
Son rôle fondateur dans la création de \textit{Fundamenta Mathematicae},
sa réflexion sur le réalisme et l’idéalisme en mathématiques, ainsi que
son engagement éducatif à travers le \textit{Poradnik dla samouków},
témoignent d’une conception unifiée du savoir, à la fois rigoureuse et
civique.
La mort prématurée de Janiszewski priva la science polonaise d’une figure
majeure, mais son héritage intellectuel survécut dans l’école polonaise
de mathématiques, qui allait rayonner internationalement durant l’entre-
deux-guerres.

\printbibliography

@online{Legionów1916,
	ALTauthor = {Legionów Polskich},
	ALTeditor = {Bureau central des publications N.K.N},
	title = {Album des Légions polonaises 126 images},
	date = {1916},
	url = {https://www.wbc.poznan.pl/dlibra/publication/162419/edition/165983/content},
	OPTsubtitle = {subtitle},
	OPTtitleaddon = {titleaddon},
	language = {Polish},
	OPTversion = {version},
	OPTnote = {note},
	organization = {Bureau central des publications N.K.N},
	OPTdate = {date},
	OPTmonth = {month},
	OPTyear = {year},
	OPTaddendum = {addendum},
	OPTpubstate = {pubstate},
	OPTurldate = {urldate},
}

@mvbook{Norman1982,
	author = {Davies, Norman},
	title = {God’s Playground: A History of Poland},
	date = {1982},
	OPTeditor = {editor},
	OPTeditora = {editora},
	OPTeditorb = {editorb},
	OPTeditorc = {editorc},
	OPTtranslator = {translator},
	OPTannotator = {annotator},
	OPTcommentator = {commentator},
	OPTintroduction = {introduction},
	OPTforeword = {foreword},
	OPTafterword = {afterword},
	subtitle = {1795 to the Present},
	OPTtitleaddon = {titleaddon},
	language = {English},
	OPToriglanguage = {origlanguage},
	OPTedition = {edition},
	volumes = {2},
	OPTseries = {series},
	OPTnumber = {number},
	OPTnote = {note},
	publisher = { New York : Columbia University Press},
	location = {New York, USA},
	isbn = {978-0231128193},
	pagetotal = {772},
	OPTaddendum = {addendum},
	OPTpubstate = {pubstate},
	OPTdoi = {doi},
	OPTeprint = {eprint},
	OPTeprintclass = {eprintclass},
	OPTeprinttype = {eprinttype},
	OPTurl = {url},
	OPTurldate = {urldate},
}

@article{Woleński2022,
	author = {Jan Woleński},
	title = {Foundations of Mathematics and Mathematical Practice. The Case of Polish Mathematical School},
	journaltitle = {Studia Historiae Scientiarum},
	date = {2022-08-26},
	OPTtranslator = {translator},
	OPTannotator = {annotator},
	OPTcommentator = {commentator},
	OPTsubtitle = {subtitle},
	OPTtitleaddon = {titleaddon},
	editor = {Michał Kokowski},
	OPTeditora = {editora},
	OPTeditorb = {editorb},
	OPTeditorc = {editorc},
	OPTjournalsubtitle = {journalsubtitle},
	OPTissuetitle = {issuetitle},
	OPTissuesubtitle = {issuesubtitle},
	language = {English},
	OPToriglanguage = {origlanguage},
	OPTseries = {series},
	volume = {21},
	number = {9},
	eid = {2543-702X},
	OPTissue = {issue},
	month = {08},
	pages = {237-257},
	OPTversion = {version},
	OPTnote = {note},
	issn = {2451-3202},
	OPTaddendum = {addendum},
	OPTpubstate = {pubstate},
	doi = {doi.org/10.4467/2543702xshs.22.007.15973},
	OPTeprint = {eprint},
	OPTeprintclass = {eprintclass},
	OPTeprinttype = {eprinttype},
	OPTurl = {url},
	OPTurldate = {urldate},
}

@article{Murawski2004,
	author = {Roman Murawski},
	title = {Philosophical reflection on mathematics in Poland in the interwar period},
	journaltitle = {Annals of Pure and Applied Logic},
	date = {2004},
	OPTtranslator = {translator},
	OPTannotator = {annotator},
	OPTcommentator = {commentator},
	OPTsubtitle = {subtitle},
	OPTtitleaddon = {titleaddon},
	editor = {Z. Adamowicz and S. Artemov and D. Niwinski and E. Orlowska and A. Romanowska and J. Wolenski},
	OPTeditora = {editora},
	OPTeditorb = {editorb},
	OPTeditorc = {editorc},
	journalsubtitle = {Provinces of logic determined. Essays in the memory of Alfred Tarski. Parts IV, V and VI},
	OPTissuetitle = {issuetitle},
	OPTissuesubtitle = {issuesubtitle},
    language = {English},
	OPToriglanguage = {origlanguage},
	OPTseries = {series},
	volume = {127},
	OPTnumber = {number},
	OPTeid = {eid},
	issue = {1-3},
	month = {06},
	pages = {325-337},
	OPTversion = {version},
	OPTnote = {note},
	issn = {0168-0072},
	OPTaddendum = {addendum},
	OPTpubstate = {pubstate},
	doi = {doi.org/10.1016/j.apal.2003.11.026},
	OPTeprint = {eprint},
	OPTeprintclass = {eprintclass},
	OPTeprinttype = {eprinttype},
	OPTurl = {url},
	OPTurldate = {urldate},
}

@incollection{Ger1986,
	author = {Céline Ger vais-Francelle},
	editor = {Éditions de la Sorbonne},
	title = {1905 La première révolution russe},
	booktitle = {La grève scolaire dans le royaume de Pologne},
	date = {1986},
	OPTeditora = {editora},
	OPTeditorb = {editorb},
	OPTeditorc = {editorc},
	OPTtranslator = {translator},
	OPTannotator = {annotator},
	OPTcommentator = {commentator},
	OPTintroduction = {introduction},
	OPTforeword = {foreword},
	OPTafterword = {afterword},
	OPTsubtitle = {subtitle},
	OPTtitleaddon = {titleaddon},
	OPTmaintitle = {maintitle},
	OPTmainsubtitle = {mainsubtitle},
	OPTmaintitleaddon = {maintitleaddon},
	OPTbooksubtitle = {booksubtitle},
	OPTbooktitleaddon = {booktitleaddon},
	language = {French},
	OPToriglanguage = {origlanguage},
	OPTvolume = {volume},
	OPTpart = {part},
	OPTedition = {edition},
	OPTvolumes = {volumes},
	OPTseries = {series},
	OPTnumber = {number},
	OPTnote = {note},
	publisher = {},
	location = {Paris, France},
	isbn = {978-2-85944-114-2},
	chapter = {IV},
	pages = {261-300},
	OPTaddendum = {addendum},
	OPTpubstate = {pubstate},
	doi = { 10.4000/books.psorbonne.52221},
	OPTeprint = {eprint},
	OPTeprintclass = {eprintclass},
	OPTeprinttype = {eprinttype},
	OPTurl = {url},
	OPTurldate = {urldate},
}

@article{Bechtel2005,
author = {Delphine Bechtel},
title = {Lemberg/Lwów/Lvov/Lviv identités d'une « ville aux frontières imprécises »},
journal = {Revue Diog\`ene(Revue internationale des sciences humaines)},
year = {2005},
volume = {2},
number = {210},
pages = {73--84}
}

@inbook{Zarnowska1986,
	author = {Anna Zarnowska and Janusz Zarnowski},
	title = {La classe ouvrière du Royaume de Pologne dans la révolution de 1905-1907},
	booktitle = {1905 La première révolution russe},
	date = {1986},
	editor = { François-Xavier Coquin and Céline Gervais-Francelle},
	language = {French},
	publisher = {Éditions de la Sorbonne},
	location = {Paris, France},
	isbn = {9782859441142},
	chapter = {La révolution dans le Royaume de Pologne},
	pages = {229--245},
	doi = {https://doi.org/10.4000/books.psorbonne.52451},
}

@article{Przeniosło2011,
 author = {Małgorzata Przeniosło},
 title = {Zygmunt Janiszewski : matematyk, legionista, filantrop (1888-1920)},
 year =  {2011},
 volume =  {18},
 number =  {1/33},
 journal =  {Niepodległość i Pamięć},
 pages =  {175-192}
}

@thesis{Janiszewski1911,
	author = {Janiszewski, Zygmunt},
	title = {Sur les continus irréductibles entre deux points},
	type = {Doctorat},
	institution = {institution},
	date = {1911-06-19},
	OPTsubtitle = {subtitle},
	OPTtitleaddon = {titleaddon},
	language = {French},
	OPTnote = {note},
	location = {Paris},
	month = {06},
	OPTisbn = {isbn},
	OPTchapter = {chapter},
	OPTpages = {pages},
	pagetotal = {93},
	OPTaddendum = {addendum},
	OPTpubstate = {pubstate},
	OPTdoi = {doi},
	OPTeprint = {eprint},
	OPTeprintclass = {eprintclass},
	OPTeprinttype = {eprinttype},
	OPTurl = {url},
	OPTurldate = {urldate},
}

@mvbook{Janiszewski1915b,
	author = {Zygmunt Janiszewski},
	title = {Logistyka, in: Poradnik dla samouków. Wskazówki metodyczne dla studiujących poszczególne nauki(Logistique, à: Poradnik dla samouków. Conseils méthodologiques pour étudier les sciences individuelles)},
	date = {1915},
	editor = {A. Heflicha and St. Michalskiego},
	OPTeditora = {editora},
	OPTeditorb = {editorb},
	OPTeditorc = {editorc},
	OPTtranslator = {translator},
	OPTannotator = {annotator},
	OPTcommentator = {commentator},
	OPTintroduction = {introduction},
	OPTforeword = {foreword},
	OPTafterword = {afterword},
	OPTsubtitle = {subtitle},
	OPTtitleaddon = {titleaddon},
	language = {Polish},
	OPToriglanguage = {origlanguage},
	OPTedition = {edition},
	volumes = {1},
	OPTseries = {series},
	OPTnumber = {number},
	OPTnote = {note},
	publisher = { L. ANCZYCA 1 SPÓŁKI},
	location = {Polish},
	OPTisbn = {isbn},
	pagetotal = {449–461},
	OPTaddendum = {addendum},
	OPTpubstate = {pubstate},
	OPTdoi = {doi},
	OPTeprint = {eprint},
	OPTeprintclass = {eprintclass},
	OPTeprinttype = {eprinttype},
	OPTurl = {url},
	OPTurldate = {urldate},
}

@article{Janiszewski1916,
	author = {Zygmunt Janiszewski},
	title = { O realizmie i idealizmie w matematyce (On realism and idealism in mathematics)},
	journaltitle = {Przegl ad Filozoficzny},
	date = {1916},
	OPTtranslator = {translator},
	OPTannotator = {annotator},
	OPTcommentator = {commentator},
	OPTsubtitle = {subtitle},
	OPTtitleaddon = {titleaddon},
	OPTeditor = {editor},
	OPTeditora = {editora},
	OPTeditorb = {editorb},
	OPTeditorc = {editorc},
	OPTjournalsubtitle = {journalsubtitle},
	OPTissuetitle = {issuetitle},
	OPTissuesubtitle = {issuesubtitle},
	language = {Polish},
	OPToriglanguage = {origlanguage},
	OPTseries = {series},
	volume = {19},
	OPTnumber = {number},
	OPTeid = {eid},
	OPTissue = {issue},
	OPTmonth = {month},
	pages = {161–170},
	OPTversion = {version},
	OPTnote = {note},
	OPTissn = {issn},
	OPTaddendum = {addendum},
	OPTpubstate = {pubstate},
	OPTdoi = {doi},
	OPTeprint = {eprint},
	OPTeprintclass = {eprintclass},
	OPTeprinttype = {eprinttype},
	OPTurl = {url},
	OPTurldate = {urldate},
}

@article{Domoradzki2018,
author = {Domoradzki, Stanisław and Stawiska, Małgorzata},
year = {2018},
month = {12},
pages = {23--49},
title = {Polish mathematicians and mathematics in World War I. Part I: Galicia (Austro-Hungarian Empire)},
volume = {17},
journal = {Studia Historiae Scientiarum},
doi = {10.4467/2543702XSHS.18.003.9323}
}
\end{document}